\documentclass[reqno,11pt]{amsart}
\usepackage{amsmath,amssymb}
\usepackage{graphicx}
\usepackage[top=1in,bottom=1in,left=1in,right=1in]{geometry}
\newtheorem{thm}{Theorem}[section]
\newtheorem{cor}[thm]{Corollary}
\newtheorem{lem}[thm]{Lemma}

\theoremstyle{definition}
\newtheorem{rem}{Remark}[section]
\numberwithin{equation}{section}

\usepackage{color}

\def\R{\mathbb R}

\def\be{\begin{equation}}
\def\ee{\end{equation}}
\def\rife#1{(\ref{#1})}
\def\vep{\varepsilon}
\def\into{\int_\Omega}
\def\elle#1{L^{#1}(\Omega)}
\def\dive{{\rm div}}
\def\de{\delta}
\def\vfi{\varphi}

\begin{document}
\title[On the regularity of  the total variation minimizers]{On the regularity of  the total variation minimizers}
\maketitle


\centerline{\bf Alessio Porretta\footnote{
Universit\`a di Roma Tor Vergata, Dipartimento di Matematica,
Via della Ricerca Scientifica 1,
00133 Roma, Italia. Email: porretta@mat.uniroma2.it}
}

\begin{abstract} We prove regularity results for the unique minimizer of the total variation functional, currently used in image processing analysis since the work by  L. Rudin, S. Osher and  E. Fatemi. In particular we show that if the source term $f$ is locally, respectively globally, Lipschitz, then the solution has the same regularity with local, respectively global,  Lipschitz   norm estimated accordingly. The result is proved in any dimension and for any (regular) domain. So far we extend  a similar result proved earlier by V. Caselles, A. Chambolle and M. Novaga  for dimension $N\leq 7$ and (in case of the global regularity)  for convex domains.
\end{abstract}

\vskip3em
\section{Introduction}

Let $\Omega$ be a  bounded smooth domain in $\R^n$, $n\geq 1$. Since the celebrated work by  L. Rudin, S. Osher and  E. Fatemi (\cite{ROS}), the problem of minimization of the total variation functional
\be\label{TV}
J(u)=\into |\nabla u| \, dx +  \frac12 \into |u-f|^2\, dx
\ee
has been extensively investigated  playing a crucial role in image processing analysis (see \cite{Br3}, \cite{CCN2}). For $f\in \elle 2$, the existence of a unique minimizer $u\in \elle2\cap BV(\Omega)$ is a standard consequence of semicontinuity of the total variation norm and properties of the space $BV(\Omega)$ consisting of functions with bounded variation (see \cite{AFP}).
 
The minimizer is unique, due to strict convexity, and is the unique solution, in a suitable weak sense, of the Neumann boundary value problem
\be\label{pb}
\begin{cases}
-\dive\left(\frac{\nabla u}{|\nabla u|}\right) +  u = f  & \qquad \hbox{in $\Omega$,}\\
\frac{\partial u}{\partial \nu}=0 & \qquad \hbox{on $\partial \Omega$.}
\end{cases}
\ee
The formulation itself, as well as the existence and uniqueness of solutions of \rife{pb} in a  large generality, has been the object of  several  papers (mostly devoted to the evolution case, see e.g. \cite{ACM} and the survey \cite{ACM2}).  Global and local regularity of the unique minimizer $u$ was investigated in a series of papers by V. Caselles, A. Chambolle and M. Novaga, see \cite{CCN3}, \cite{CCN} and the survey \cite{CCN2}. Among their results, the authors proved that the solution $u$ is locally  H\"older or Lipschitz continuous whenever so is the source $f$, at least  in dimension $n\leq 7$. In addition, if $f$ is globally   H\"older or Lipschitz in $\Omega$,  the corresponding  global regularity  was also proved for $u$ assuming  that $\Omega$ is a convex domain.   Only in the very recent paper \cite{M}, some of these results were extended to the case of any dimension $n$; in particular G. Mercier has proved that the  continuity of $f$ implies the continuity of $u$ and, in the case of convex domains, the  modulus of continuity is also inherited globally by the solution.
The issue of regularity of solutions was  newly addressed  by  recent results of H. Brezis (\cite{Br}) and  T. Sznigir  (\cite{Sz}), showing further regularity of the solution in the one-dimensional case. In particular,    whenever $n=1$ and $\Omega$ is an interval, it is proved in \cite{Br} that, at least for smooth $f$, one has  $Du\in BV(\Omega)$.

\vskip1em
The purpose of this article is to extend some results proved in  \cite{CCN} and to introduce  a new possible strategy for proving regularity of solutions. The main contribution that we give is to show that $f$ locally Lipschitz implies $u$ locally Lipschitz in {\it any} dimension; additionally, if $f$ is globally Lipschitz, then $u$ is globally Lipschitz in {\it any} domain. In particular, our results extend  the local or, respectively,  global Lipschitz estimate obtained in \cite{CCN}  to any dimension and to any domain. Our approach also shows why the convexity of the domain plays a role and indeed we prove further results (Sobolev regularity estimates for the solution) but only in convex domains.

\vskip1em
Let us now state our  main results. We start by the preservation of the local Lipschitz regularity.  

\begin{thm}\label{liploc} Let $u\in BV(\Omega)\cap\elle2$ be the minimizer of \rife{TV}. If $f$ is locally Lipschitz in $\Omega$, then  $u$ is locally Lipschitz in $\Omega$ and  there exists $C>0$ such that 
$$
\|\nabla u\|_{L^\infty(B_R(x_0))} \leq   \left( \|\nabla f\|_{L^\infty(B_{2R}(x_0))} +   \frac C{R^2}\right)
$$
for any $x_0\in  \Omega$, $R>0$ such that $B_{2R}(x_0)\subset \Omega$ and  $f$ is Lipschitz in $B_{2R}(x_0)$.
\end{thm}

\vskip1em
Global Lipschitz continuity of $f$ yields global Lipschitz regularity of $u$ as well. Indeed we have:

\begin{thm}\label{lip} Let $u\in BV(\Omega)\cap\elle2$ be the minimizer of \rife{TV}. If $f\in W^{1,\infty}(\Omega)$, then we have $u\in W^{1,\infty}(\Omega)$ and 
$$
\|\nabla u\|_{L^\infty(\Omega)} \leq  c_1\, \left( \|\nabla f\|_{L^\infty(\Omega)}+ c_0\right)\,,  
$$
where $c_0,c_1$ only depend  on $\Omega$. 
\end{thm}
\vskip1em
A special case is given by convex domains, where the above result holds with $c_0=0$ and $c_1=1$, namely there is preservation of the Lipschitz norm in this case. We prove that this is actually true for all  norms in Sobolev spaces $W^{1,p}(\Omega)$ with $p\geq 2$.

\begin{thm}\label{sob} Let $u\in BV(\Omega)\cap\elle2$ be the minimizer of \rife{TV}, and assume that $\Omega$ is convex. For every $p\in [2,+\infty]$, we have that if $f\in W^{1,p}(\Omega)$ then $u\in W^{1,p}(\Omega)$ and 
$$
\|\nabla u\|_{L^p(\Omega)} \leq    \|\nabla f\|_{L^p(\Omega)}\,.  
$$
\end{thm}

Unfortunately, we are unable to extend the above result to $p<2$.  We point out that Theorem \ref{sob} holds in dimension one for any value of $p>1$, as well as for the BV norm, as a consequence of the representation formula provided in \cite{Br} for the minimizers, see Remark \ref{d=1}.

\vskip0.5em
Let us stress that {\it all the estimates obtained  in Theorems \ref{liploc}-\ref{sob} are scale invariant and independent of the coercivity of  the second order term}, often considered as the regularizing term in image processing (in the so-called denoising problem). Namely, if we  consider the one-parameter family of  functionals
\be\label{TVmu}
J_\mu(u):= \mu\, \into |D u| \, dx +  \frac1 2 \into |u-f|^2\, dx  
\ee
then all the above estimates still hold in the same form for minimizers of $J_\mu$, namely the estimates hold independently of $\mu$ (say, for all $\mu$ in a  bounded set). This is easily observed by scaling properties. 
\vskip0.5em

Let us now comment about the methods used so far for proving regularity of the total variation minimizer.
On one hand,  the approach used in \cite{Br}, \cite{Sz} for the  one-dimensional  results relies on  a duality argument which seems not applicable to problem \rife{pb} in higher dimension. On another hand, the method used by  V. Caselles, A. Chambolle and M. Novaga relies on the characterization of the level sets of the solution as minimal sets for  the prescribed curvature problem. On account of this characterization, their regularity results stand on previously established regularity  of optimal sets for the prescribed curvature problem; those results are responsible for the limitation $n\leq 7$ which is required in their approach.

By contrast, the strategy that we use seems very natural and applies to any dimension; it dates back to a classical method for gradient estimates, often referred to as Bernstein's method (\cite{ber}). This tool, refined by several authors later, (see e.g. \cite{GT}, \cite{L1}, \cite{Se}), proved successful for quasilinear equations in both divergence and non-divergence form. It is not by  chance that a similar approach was used for the prescribed  or the mean curvature problem in different contexts, see e.g. \cite{Ba}, \cite{Wa}, as well as for $p$-Laplace type equations (\cite{LP2}).  

This method relies on the fact that, if $u$ is a solution to a  quasilinear elliptic equation, then  $|\nabla u|^2$ can be proved to be a  sub solution to an elliptic equation, so that it can  be estimated through maximum principle. 
We show that this  idea can be exploited for the total variation minimizer even if the Euler equation  is strongly degenerate. Indeed, we show that  $|\nabla u|^2$ can be regarded both as a  sub solution of a  non-divergence form equation (which yields estimates through maximum principle) and as  a  sub solution of a divergence form equation (which yields integral estimates). However, and this explains the difference between the above statements,  while the estimate on the maximum norm can be  localized, this is not so (clear) for the integral estimates. This explains the stronger conditions required for the estimate of Sobolev norms (Theorem \ref{sob}); in fact, our proof readily shows how the convexity of the domain is   a  very natural condition which allows for global estimates avoiding the issue of localization.

Of course, in order to carefully apply this kind of arguments,  we will   prove uniform a priori estimates for the unique smooth solution 
%
%
$u_{\vep,\de}$   of the regularized problem
$$
\begin{cases}
-\de \Delta u  -\dive\left(\frac{\nabla u}{\sqrt{\vep +|\nabla u|^2}}\right) +   u = f  & \qquad \hbox{in $\Omega$,}\\
\frac{\partial u}{\partial \nu}=0 & \qquad \hbox{on $\partial \Omega$.}
\end{cases}
$$
We notice that the role of the two parameters $\vep, \de$ is somehow different: $\de$ is used in order to recover a uniformly elliptic operator and to perform our estimates using smooth solutions. The parameter $\vep$ is rather used   in order to avoid the singularity of the energy functional.  In all our estimates, but for the local Lipschitz regularity, there would be no loss of generality in taking $\de=\vep$, however it will be more convenient to carry on the analysis with possibly different parameters. Eventually, with rather simple stability properties, we will  deduce from such approximations the  estimates for the unique minimum of \rife{TV}.

As a final comment, it is worth  pointing out that  there exist generalized versions of the classical Bernstein's method which replace the gradient estimate with oscillation estimates and can therefore be handled in a less regular framework, just  through comparison principles. This might possibly establish a connection between the purely PDE approach we present here and the more geometrical one proposed in \cite{M}  as a refinement of the slicing method used by V. Caselles, A. Chambolle and M. Novaga.



\section{Local and global Lipschitz estimates}

In this Section we derive the Lipschitz estimates for the regularized equation
\be\label{pbeps}
\begin{cases}
-\de \Delta u  -\dive\left(\frac{\nabla u}{\sqrt{\vep +|\nabla u|^2}}\right) +  \lambda\,  u = f  & \qquad \hbox{in $\Omega$,}\\
\frac{\partial u}{\partial \nu}=0 & \qquad \hbox{on $\partial \Omega$.}
\end{cases}
\ee
Even if the parameter $\lambda$ could be rescaled, it will be useful to distinguish the zero order term from first and second order terms in the computations below, where we establish the equation satisfied by the squared gradient of the solution.
We point out that, due to the elliptic regularization,  solutions of \rife{pbeps} are smooth enough to justify the following computations, see \cite{GT}.  

Even if the solution depends on $\vep$, $\de$, we will eventually drop the index and write $u$ for $u_{\vep, \de}$.  We denote $|D^2u|^2= \sum_{i,j} u_{ij}^2$, where $u_{ij}= \frac{\partial^2 u}{\partial x_i\partial x_j}$. Henceforth, we use $C$ to denote possibly different real numbers which are independent of the parameters $\vep,\de$ as well as of $u$, $\lambda$  and $f$. The constant $C$ may sometimes depend on $\Omega$, but it will be made clear in the proofs if this is the case.

\begin{lem}\label{base}
Assume that $f\in W^{1,\infty}(\Omega)$. Let $u=u_{\vep,\de}$ be the unique solution of \rife{pbeps}. Then, the function $w:= |\nabla u_{\vep,\de}|^2$ satisfies the equation
\be\label{eqw}
\begin{split}
{\mathcal L}  w & + 2\lambda\, w+ 2\de \, |D^2u|^2 +
2\frac{|D^2 u|^2}{{\sqrt{\vep +w}}} 
\\
& = -\frac{\Delta u}{{(\vep +w)^{\frac32}}}\, \left(\nabla w \nabla u\right)
+\frac32 \frac{\left(\nabla u \nabla  w\right)^2}{(\vep+ w)^{\frac52}}
-\frac12\frac{|\nabla  w|^2}{(\vep+ w)^{\frac32}}+ 2\nabla f\, \nabla u
\end{split}
\ee
where the operator ${\mathcal L}$ is defined as
$$
 {\mathcal L} w:  = -\de \Delta w - \frac{\Delta w }{{\sqrt{\vep +|\nabla u|^2}}} + \frac{D^2w \nabla u\nabla u}{(\vep+ |\nabla u|^2)^{\frac32}}\,.
$$
As a consequence, $w$ satisfies the following inequality
\be\label{redu}
{\mathcal L}  w   + 2\lambda\, w \leq C \frac{|\nabla  w|^2}{(\vep+ w)^{\frac32}} + 2|\nabla f|\, \sqrt w\,.
\ee
%
\end{lem}

\begin{rem}
Notice that the operator ${\mathcal L}$ is elliptic. In particular, the lemma says that $w$ is  a positive sub solution of an elliptic equation.  
\end{rem}

\proof
Let us first notice that \rife{pbeps} can be read as
\be\label{svilup}
\de \, \Delta u + \frac{\Delta u }{{\sqrt{\vep +w}}} - \frac12 \frac{\nabla u \nabla  w}{(\vep+ w)^{\frac32}}= \lambda u -f \,.
\ee
Moreover, we have:
\be\label{one}
\Delta w = 2 \nabla u \nabla (\Delta u) +2 |D^2u |^2\,.
\ee
From \rife{svilup}, we have
\begin{align*} &
\de \, \nabla (\Delta u)\nabla u + \frac{\nabla (\Delta u) \nabla u }{{\sqrt{\vep +w}}}  =
\frac12\frac{\Delta u}{{(\vep +w)^{\frac32}}}\, \left(\nabla w \nabla u\right) 
\\ &  \quad 
-\frac34 \frac{\left(\nabla u \nabla  w\right)^2}{(\vep+ w)^{\frac52}}
+\frac12 \frac{D^2 w\, \nabla u\nabla u}{(\vep+ w)^{\frac32}}
+\frac14\frac{|\nabla  w|^2}{(\vep+ w)^{\frac32}}
\\
& \quad + \lambda |\nabla u|^2- \nabla f\nabla u
\end{align*}
where we used that $D^2u \nabla u= \frac12 \nabla w$.  Using equality \rife{one} in the first two terms, we obtain
\begin{align*} &
\frac\de2 \, \Delta w + \frac12\frac{\Delta w }{{\sqrt{\vep +w}}}  = \de \, |D^2(u)|^2 +
\frac{|D^2 u|^2}{{\sqrt{\vep +w}}}
+ \frac12\frac{\Delta u}{{(\vep +w)^{\frac32}}}\, \left(\nabla w \nabla u\right) 
\\ &  \quad 
-\frac34 \frac{\left(\nabla u \nabla  w\right)^2}{(\vep+ w)^{\frac52}}
+\frac12 \frac{D^2 w\, \nabla u\nabla u}{(\vep+ w)^{\frac32}}
+\frac14\frac{|\nabla  w|^2}{(\vep+ w)^{\frac32}}
\\
& \quad + \lambda |\nabla u|^2- \nabla f\nabla u\,.
\end{align*}
Now, if we  define the  operator
\begin{align*}
{\mathcal L}  w: & = -\de \Delta w - \frac{\Delta w }{{\sqrt{\vep +|Du|^2}}} + \frac{D^2w \nabla u\nabla u}{(\vep+ |Du|^2)^{\frac32}}
\\
& = 
-\de \Delta w - \frac{\Delta w }{{\sqrt{\vep +w}}} + \frac{D^2w \nabla u\nabla u}{(\vep+ w)^{\frac32}}\,,
\end{align*}
we have derived the following equation satisfied by $w$:
\begin{align*}
{\mathcal L}  w & + 2\lambda\, w+ 2\de \, |D^2u|^2 +
2\frac{|D^2 u|^2}{{\sqrt{\vep +w}}} 
\\
& = -\frac{\Delta u}{{(\vep +w)^{\frac32}}}\, \left(\nabla w \nabla u\right)
+\frac32 \frac{\left(\nabla u \nabla  w\right)^2}{(\vep+ w)^{\frac52}}
-\frac12\frac{|\nabla  w|^2}{(\vep+ w)^{\frac32}}+ 2\nabla f\, \nabla u
\end{align*}
which is \rife{eqw}.
%

Now we further estimate  the right-hand side. Indeed, using Young's inequality and recalling that $|\nabla u|^2=w$, we have
\begin{align*}  
-\frac{\Delta u}{{(\vep +w)^{\frac32}}}\, \left(\nabla w \nabla u\right)
+\frac32 \frac{\left(\nabla u \nabla  w\right)^2}{(\vep+ w)^{\frac52}}
& \leq  \frac12\frac{|D^2u|^2}{{(\vep +w)^{\frac12}}}+ C \frac{\left(\nabla u \nabla  w\right)^2}{(\vep+ w)^{\frac52}}
\\
& \leq  \frac12\frac{|D^2u|^2}{{(\vep +w)^{\frac12}}}+ C \frac{|\nabla  w|\, |D^2 u|\,  |\nabla u|}{(\vep+ w)^{\frac32}} 
\\
& \leq
 \frac{|D^2u|^2}{{(\vep +w)^{\frac12}}}+ C \frac{|\nabla  w|^2}{(\vep+ w)^{\frac32}}\,.
\end{align*}
Therefore, we deduce that
$$
{\mathcal L}  w   + 2\lambda\, w + \frac{|D^2u|^2}{{(\vep +w)^{\frac12}}} \leq C \frac{|\nabla  w|^2}{(\vep+ w)^{\frac32}} + 2|\nabla f|\, |\nabla u|\, 
$$
which implies \rife{redu}.
\qed

\vskip1em

Lemma \ref{base} says that $w$ is  a positive sub solution of an elliptic equation, and so $w$ obeys to the maximum principle. In particular, we could already deduce a global bound in convex domains, since the convexity of $\Omega$ and the Neumann condition $\frac{\partial u}{\partial \nu}=0$ imply that $\frac{\partial w}{\partial \nu}\leq 0$ (see Lemma \ref{neuma} below). 

In order to get at local results, as well as to drop the convexity condition on the domain, 
we need a  possibly localized version of the  inequality \rife{redu}.

\begin{lem}\label{loc} Assume that $f\in W^{1,\infty}(\Omega)$. Let $u_{\vep,\de}$ be the unique solution of \rife{pbeps}, and define as before $w:= |\nabla u_{\vep,\de}|^2$. Then, for any $\vfi\in C^2(\Omega)$ the following inequality
holds in the open set $\{x\in \Omega:\,\, \vfi(x)>0\}$:
\be\label{ineloc}
\begin{split}
& {\mathcal L}(w\vfi)+ 2\lambda\, w \,\vfi   \leq C\, \frac{|\nabla  (w\vfi)|^2}{\vfi(\vep+ w)^{\frac32}}
- 2 \de \nabla (w\vfi)\frac{\nabla \vfi}{\vfi}+   \de \, w \left[\frac{2|\nabla \vfi|^2}{\vfi}-\Delta \vfi \right]
\\
& \quad +  2|\nabla f|\, \sqrt w\, \vfi 
+  C\, \sqrt w\, \left[\frac{|\nabla \vfi|^2}{\vfi} + |D^2 \vfi| \right]\, .
\end{split}
\ee
\end{lem}

\proof Let $\vfi$ be a $C^2$ function. Since we have
$$
(w\vfi)_{ij}= w_{ij}\, \vfi +  (w_i\, \vfi_j+ w_j\, \vfi_i )  + w \vfi_{ij}
$$
then
\begin{align*}
\left({\mathcal L}  w\right)\, \vfi & = {\mathcal L}(w\vfi)+ 2\de \nabla w\nabla \vfi + \de \, w\, \Delta \vfi + 2  \frac{\nabla w\nabla \vfi  }{{\sqrt{\vep +|\nabla u|^2}}} + \frac {w\, \Delta \vfi}{{\sqrt{\vep +|\nabla u|^2}}}  
\\
& 
\quad 
 - 2\frac{(\nabla w\nabla u) (\nabla \vfi\nabla u)}{(\vep+ |\nabla u|^2)^{\frac32}}
- \frac w{(\vep+ |\nabla u|^2)^{\frac12}}\, \frac{D^2\vfi \nabla u\nabla u}{(\vep+ |\nabla u|^2)}\,.
\end{align*}
Then from \rife{redu} we obtain
\begin{align*}
& {\mathcal L}(w\vfi)+ 2\lambda\, w \,\vfi \leq C \frac{|\nabla  w|^2\, \vfi}{(\vep+ w)^{\frac32}} + 2|\nabla f|\, \sqrt w\, \vfi
\\
& \quad - 2\de \nabla w\nabla \vfi - \de \, w\, \Delta \vfi - 2  \frac{\nabla w\nabla \vfi  }{{\sqrt{\vep +|\nabla u|^2}}}  - \frac {w\, \Delta \vfi}{{\sqrt{\vep +|\nabla u|^2}}}  
\\
& 
\quad  + 2\frac{(\nabla w\nabla u) (\nabla \vfi\nabla u)}{(\vep+ |\nabla u|^2)^{\frac32}}
+ \frac w{(\vep+ |\nabla u|^2)^{\frac12}}\, \frac{D^2\vfi \nabla u\nabla u}{(\vep+ |\nabla u|^2)}\,.
\end{align*}
We estimate further the terms in the second and third line; by Young's inequality, and using  that $|\nabla u|= \sqrt w$ we get
\begin{align*}
& {\mathcal L}(w\vfi)+ 2\lambda\, w \,\vfi \leq C \frac{|\nabla  w|^2\, \vfi}{(\vep+ w)^{\frac32}} + 2|\nabla f|\, \sqrt w\, \vfi
\\
& \quad - 2\de \nabla w\nabla \vfi - \de \, w\, \Delta \vfi 
+  C\, \sqrt w\, \left[\frac{|\nabla \vfi|^2}{\vfi} + |D^2 \vfi| \right]\,.
\end{align*}
Now we use that 
$$
2\de \nabla w\nabla \vfi = 2 \de \nabla (w\vfi)\frac{\nabla \vfi}{\vfi}- 2\de \, w\frac{|\nabla \vfi|^2}{\vfi}
$$
and similarly
$$
\frac{|\nabla  w|^2\, \vfi}{(\vep+ w)^{\frac32}} \leq 2\frac{|\nabla  (w\vfi)|^2}{\vfi(\vep+ w)^{\frac32}} + 2\, \sqrt w\, \frac{|\nabla    \vfi|^2}{\vfi}\,. 
$$
Putting all  together   we get
\begin{align*}
& {\mathcal L}(w\vfi)+ 2\lambda\, w \,\vfi   \leq C\, \frac{|\nabla  (w\vfi)|^2}{\vfi(\vep+ w)^{\frac32}}
- 2 \de \nabla (w\vfi)\frac{\nabla \vfi}{\vfi} + \de\, w \left[ \frac{2|\nabla \vfi|^2}{\vfi}-\Delta \vfi\right]
\\
& \quad +  2|\nabla f|\, \sqrt w\, \vfi 
+  C\, \sqrt w\, \left[\frac{|\nabla \vfi|^2}{\vfi} + |D^2 \vfi| \right]\,.
\end{align*}
\qed

\vskip1em
Now we use the Neumann boundary condition to get at global estimates. The following useful lemma is classical, see e.g. \cite{L2}, \cite{LP2}. For the reader's convenience, we give shortly the proof below.  We denote by $d(x)$   the distance function to the boundary; since $\Omega$ is smooth we have that  $d(x) $ is a $C^2$ function in a  neighborhood of the boundary and $\nabla d(x)= - \nu(x)$ on $\partial \Omega$, where $\nu$ is the outward unit normal.

\begin{lem}\label{neuma} Let $u_{\vep,\de}$ be the unique solution of \rife{pbeps}, and define as before $w:= |\nabla u_{\vep,\de}|^2$. 

If $\Omega$ is convex, then $w$ satisfies $\frac{\partial w}{\partial \nu}  \leq 0$ on $\partial \Omega$. 
More generally, if $z= w\, e^{\gamma\, d(x)}$, we have $\frac{\partial z}{\partial\nu}  \leq 0$ on $\partial \Omega$
for any $\gamma\geq 2\|(D^2d)_+\|_{L^\infty(\partial\Omega)}$.
\end{lem}

\proof We look at the   boundary condition for the function $w$. As before, we drop the indexes. We have
\begin{align*}
\frac{\partial w}{\partial \nu}=\nabla |\nabla u|^2 \cdot \nu &  = 2 D^2 u \nabla u\cdot \nu 
\\
& = 2 \nabla u \nabla (\nabla u\cdot \nu)- 
2 D\nu \nabla u \nabla u\,.
\end{align*}
Since $\nabla u\cdot \nu$ vanishes on the boundary, its gradient points in the normal direction, so 
$$
\nabla u \nabla (\nabla u\cdot \nu)=0
$$
using once more the Neumann condition for $u$. Therefore, we conclude that 
$$
\nabla |\nabla u|^2 \cdot \nu = - 2 D\nu \nabla u \nabla u\,.
$$
Since $\Omega$ is a smooth set, we have $D\nu(x)= - D^2 d(x)$, where $d(x)$ is the distance function to the boundary. 
First we notice that if $\Omega$ is convex then $d(x)$ is concave and $D^2 d(x)\leq 0$ in a  neighborhood of the boundary. In this case we deduce that 
$$
\frac{\partial w}{\partial \nu} \leq 0 \qquad \hbox{on $\partial \Omega$.}
$$
For   general $\Omega$, we call $z= w\, e^{\gamma\, d(x)}$. Since we have
$$
\frac{\partial z}{\partial \nu} =   e^{\gamma\, d(x)}\left( 2 D^2d(x) \nabla u \nabla u - \gamma\, |\nabla u|^2\right)\qquad x\in \partial \Omega\,,
$$
if we choose $\gamma\geq  2\|(D^2d)_+\|_{L^\infty(\partial\Omega)}$ then we have  $\frac{\partial z}{\partial\nu} \leq 0$.
\qed

\vskip1em
We deduce now  the following global estimate.

\begin{cor}\label{glob-lip} Assume that $f\in W^{1,\infty}(\Omega)$. Let $u_{\vep,\de}$ be the unique solution of \rife{pbeps}. Then there exist   constants $c_0, c_1$, only depending on $\Omega$, such that, for all $\vep$ and $\de$ sufficiently small we have  
$$
\| \nabla u_{\vep,\de}\|_\infty \leq \frac {c_1}\lambda\, \left(\|\nabla f\|_\infty+ c_0\right)\,.
$$
In addition, if $\Omega$ is convex, the inequality is true with $c_0=0$ and $ c_1=1$.
\end{cor}

\proof
In the following, we  still call $d(x)$  a $C^2$ function which is positive in $\Omega$ and coincides with the distance function in a  neighborhood of the boundary. Let us set $z= w\, e^{\gamma\, d(x)}$, with $\gamma \geq\|(D^2d)_+\|_{L^\infty(\partial\Omega)}$.  
Now we use Lemma \ref{loc} with $\vfi= e^{\gamma\, d(x)}$, and we   look at the elliptic equation satisfied by $z$ given by Lemma \ref{loc}.  Since  we easily estimate
$$
\left[\frac{|\nabla \vfi|^2}{\vfi} + |D^2 \vfi| \right] \leq C\, (\gamma^2+\gamma )\,   e^{\gamma\, d(x)} 
$$
for some constant $C$ only depending on $d(x)$, from \rife{ineloc} we deduce that $z$ satisfies in $\Omega$ the elliptic inequality 
\begin{align*}
  {\mathcal L}z+ 2\lambda\, z &  \leq C\, e^{-\gamma d}\frac{|\nabla z|^2}{(\vep+ w)^{\frac32}}
- 2 \de\, \gamma\,  \nabla z \nabla d+   
\\
& \quad +  2|\nabla f|\, \sqrt z\, e^{\frac12\gamma d} +  C\, (\gamma^2+\gamma ) (\sqrt w\, e^{\gamma d}+\de z) \,.
\end{align*}
Due to Hopf lemma, and on account of  Lemma \ref{neuma},  $z$ cannot assume its maximum at the boundary.
Therefore, on the maximum point of $z$ we have
$$
2\lambda z \leq   \sqrt z\, e^{\frac12\gamma d} \left[2|\nabla f| +  C\, (\gamma^2+\gamma )\right] + \de\, C\, (\gamma^2+\gamma ) z .
$$
For a  sufficiently small $\de$ (only depending on $\Omega$), we deduce that 
$$
\max \, \sqrt z \leq \frac {e^{\frac\gamma 2\|d\|_\infty }}{\lambda-\de\, C  (\gamma^2+\gamma )} \left(\|\nabla f\|_{\infty} + C(\gamma+\gamma^2)\right)\,,
$$
for some $C$ only depending on $\Omega$.  Since $z\geq |\nabla u|^2$, this implies the conclusion. Moreover,  if $\Omega$ is convex then  we can choose $\gamma=0$ (i.e. $\vfi=1$). Then we get at the same conclusion obtaining now
$$
\| \nabla u\|_\infty \leq \frac{\|\nabla f\|_\infty}\lambda \,.
$$
\qed

Now we prove the  local Lipschitz estimates in the limit as $\de\to0$. We notice that, if $f\in W^{1,\infty}(\Omega)$, Corollary 
\ref{glob-lip} implies that $u_{\vep, \de}$ is bounded in $W^{1,\infty}(\Omega)$, and then relatively compact in the uniform topology.

\begin{cor}\label{est-loc} Assume that $f\in W^{1,\infty}(\Omega)$.  Let $u_\vep$ be the limit of $u_{\vep, \de}$ for some subsequence $\de\to 0$.  Then, for every $x_0\in \Omega$, every  $R,\rho>0$ with  $(1+\rho)R<{\rm dist}(x_0,\partial\Omega)$,  we have 
$$
\sup_{B_R(x_0)} \, |\nabla u_\vep| \leq \frac1\lambda  \left( \frac {C_\rho}{R^2} + \sup_{B_{(1+\rho)R}(x_0)} \, |\nabla f|\right)\,,
$$
for some constant $C_\rho$ depending on $\rho$ but  independent of $\vep$ and $R$. 
\end{cor}

\proof We consider the case $R=1$, the general statement be recovered by scaling.  Let us take a cut-off function $\vfi\in C^2$ such that $\vfi=1$ in $B_1$, $\vfi$ has support in $B_{1+\rho}$ and satisfies
$$
|\nabla \vfi|^2\leq C \, \vfi^{\frac32}\,,\qquad  |D^2\vfi| \leq C \, \vfi^{\frac12}\,,
$$
for some constant $C$ depending on $\rho$. 
If $u_{\vep,\de}$ denotes the solution of problem \rife{pbeps}, we use Lemma \ref{loc} with such a  $\vfi$ and we look at the positive maximum attained by $z:= w\,\vfi$ in the ball $B_{1+\rho}$, where $w= |\nabla u_{\vep,\de}|^2$. 
On account of \rife{ineloc}, at the maximum point we get
$$
2\lambda\, z  \leq C\, \de \, w \sqrt \vfi + \sqrt z\, [  2|\nabla f|\, \sqrt \vfi +  C]\,,
$$ 
hence
$$
2\lambda\, \sqrt z  \leq C\, \de \, \|w\|_\infty^{\frac12} +    2 |\nabla f| +  C \,.
$$
This implies that 
$$
\sup_{B_1} |Du_{\vep, \de}|\leq \max_{ {B_{1+\rho}}} \sqrt z \leq \frac 1{\lambda} \left[\de\, C\,  \| Du_{\vep,\de}\|_\infty+ \sup_{B_{1+\rho}} |\nabla f| + C\right]\,.
$$
Since $f \in W^{1,\infty}(\Omega)$,   we know that $\| Du_{\vep,\de}\|_\infty$ is uniformly bounded (from Corollary \ref{glob-lip}). Hence, $u_{\vep,\de}$  is  relatively compact in the uniform topology, and $Du_{\vep,\de}$ in the weak-$*$ topology.
 For any $u_\vep$ obtained as limit of a  subsequence of $u_{\vep,\de}$, we get as $\de\to 0$: 
$$
\sup_{B_1} |Du_{\vep}|\leq   \frac 1{\lambda} \left[ \sup_{B_{1+\rho}} |\nabla f| + C\right]\,.
$$
\qed

\begin{rem} We observe that if we consider the $\vep$- approximation of the functional $J_\mu$ defined in \rife{TVmu}, then the estimate of Corollary \ref{est-loc} will be rescaled into the following one:
$$
\sup_{B_R(x_0)} \, |\nabla u_\vep^\mu| \leq \frac1\lambda  \left( \mu\frac {C_\rho}{R^2} + \sup_{B_{(1+\rho)R}(x_0)} \, |\nabla f|\right)\,,
$$
for the corresponding minimizers $u_\vep^\mu$. In particular, one can notice that, as the regularizing parameter $\mu$ vanishes, the estimate will be precised  into $\sup_{B_R(x_0)} \, |\nabla u_\vep| \leq \frac1\lambda \, \sup_{B_{(1+\rho)R}(x_0)} \, |\nabla f|$. 
By letting $\rho\to 0$ and averaging in $R$, this eventually leads  to the expected point wise estimate 
$|\nabla u| \leq \frac1\lambda \,  |\nabla f| $.
\end{rem}

\section{Integral estimates for the gradient.}

In this Section we further investigate the equation satisfied by the squared gradient, which is here observed in divergence form. The conservative character of this equation may be crucial in obtaining integral estimates for the gradient and preservation of  Sobolev norms.

\begin{lem} Let $u_{\vep,\de}$ be the unique solution of \rife{pbeps}, we define as before $w:= |\nabla u_{\vep,\de}|^2$. Then $w$ satisfies
\be\label{div-form}
\begin{split}
& -  \dive\left(\de \nabla w + \frac{\nabla w}{(\vep +w)^{\frac12}} - \frac{\left(\nabla w \nabla u\right)\nabla u}{(\vep +w)^{\frac32}}\right)  + 2\lambda\, w
\\ & \quad 
+ 2\de \, |D^2u|^2 +
2\frac{|D^2 u|^2}{{\sqrt{\vep +w}}} 
= \frac12\frac{|\nabla  w|^2}{(\vep+ w)^{\frac32}}+ 2\nabla f\, \nabla u\,.
\end{split}
\ee
\end{lem}

\proof  Let ${\mathcal L} $ be the operator defined in Lemma \ref{base}. Then we have
\begin{align*} 
& 
{\mathcal L}  w= -\de \Delta w - \frac{\Delta w }{{\sqrt{\vep +w}}} + \frac{D^2w \nabla u\nabla u}{(\vep+ w)^{\frac32}}
 = 
-  \dive\left(\de \nabla w + \frac{\nabla w}{(\vep +w)^{\frac12}} - \frac{\left(\nabla w \nabla u\right)\nabla u}{(\vep +w)^{\frac32}}\right)  
\\
& \quad - \frac12\frac{|\nabla w|^2}{(\vep +w)^{\frac32}} - 
\frac{\left(\nabla w \nabla u\right)\Delta u}{(\vep +w)^{\frac32}} 
- \frac{ D^2 u \, \nabla w \cdot \nabla u}{(\vep +w)^{\frac32}}
+ \frac32 \frac{\left(\nabla w \nabla u\right)^2}{(\vep +w)^{\frac52}} \,.
\end{align*}
Since $D^2 u \, \nabla w \cdot \nabla u= \frac12 |\nabla w|^2$, we deduce
\begin{align*} 
& 
{\mathcal L}  w
 = 
-  \dive\left(\de \nabla w + \frac{\nabla w}{(\vep +w)^{\frac12}} - \frac{\left(\nabla w \nabla u\right)\nabla u}{(\vep +w)^{\frac32}}\right)  
\\
& \quad -  \frac{|\nabla w|^2}{(\vep +w)^{\frac32}} - 
\frac{\left(\nabla w \nabla u\right)\Delta u}{(\vep +w)^{\frac32}} 
+ \frac32 \frac{\left(\nabla w \nabla u\right)^2}{(\vep +w)^{\frac52}} \,.
\end{align*}
As a consequence, from \rife{eqw} we obtain that $w$ satisfies
\begin{align*}
& -  \dive\left(\de \nabla w + \frac{\nabla w}{(\vep +w)^{\frac12}} - \frac{\left(\nabla w \nabla u\right)\nabla u}{(\vep +w)^{\frac32}}\right)  + 2\lambda\, w\\
& \quad 
+ 2\de \, |D^2u|^2 +
2\frac{|D^2 u|^2}{{\sqrt{\vep +w}}} 
= \frac12\frac{|\nabla  w|^2}{(\vep+ w)^{\frac32}}+ 2\nabla f\, \nabla u
\end{align*}
which is \rife{div-form}.
\qed

\vskip1em
We immediately deduce the following global estimate in convex domains.

\begin{cor}\label{estint}
Let $\Omega$ be a convex domain, and let $u_{\vep,\de}$ be the unique solution of \rife{pbeps}. Assume that $f\in W^{1,\infty}(\Omega)$. Then, for any $p\in [2, \infty]$ we have
\be\label{pres}
\| \nabla u_{\vep,\de}\|_{\elle p} \leq \frac 1\lambda \, \| \nabla f\|_{\elle p}\,.
\ee 
\end{cor}

\proof  As usual, we drop the indexes during the proof and we set $w:= |\nabla u|^2$. 
We first notice that, using $\nabla w= 2D^2 u \nabla u$, we have 
$$
\frac12\frac{|\nabla  w|^2}{(\vep+ w)^{\frac32}} \leq 2 \frac{|D^2u|^2\, |\nabla u|^2}{(\vep+ w)^{\frac32}} \leq 2 \frac{|D^2u|^2}{(\vep+ w)^{\frac12}}\,.
$$
Therefore, \rife{div-form} implies
\be\label{new}
-  \dive\left(\de \nabla w + \frac{\nabla w}{(\vep +w)^{\frac12}} - \frac{\left(\nabla w \nabla u\right)\nabla u}{(\vep +w)^{\frac32}}\right)  + 2\lambda\, w \leq  2|\nabla f|\, \sqrt w\,.
\ee
Now we recall from Lemma \ref{neuma} that  we have $\frac{\partial w}{\partial \nu}\leq 0$ on the boundary. In addition, the Neumann condition holds for $u$. So, when we multiply inequality \rife{new} by $w^\beta$, for $\beta\geq 0$, we deduce
$$
\beta \into \left(\frac{|\nabla w|^2}{(\vep +w)^{\frac12}} - \frac{\left(\nabla w \nabla u\right)^2}{(\vep +w)^{\frac32}}\right) w^{\beta-1}\, dx +  2\lambda \into w^{\beta+1} \leq 2 \into |\nabla f|\, w^{\beta+\frac12}\, dx\,. 
$$
The first integral is nonnegative.  Hence we get
$$
\lambda \into w^{\beta+1} \leq  \| \nabla f\|_{\elle{2(\beta+1)}}\, \left(\into w^{\beta+1}\, dx\right)^{1-\frac1{2(\beta+1)}}\, 
$$
which yields
$$
\lambda \left(\into w^{\beta+1}\right)^{\frac1{2(\beta+1)}} \leq  \| \nabla f\|_{\elle{2(\beta+1)}}\,.
$$
Recalling that $w= |\nabla u|^2$, we conclude. 
\qed

\vskip1em

\begin{rem}\label{d=1}
At this stage, we are unable to prove  a similar estimate as \rife{pres}   for   Sobolev norms with $p<2$. This is known to be true in dimension one (\cite{Br}); in that case,  H. Brezis provides with a representation formula for the minimizer, roughly speaking one has $u= v'+ f$ where $v$ solves  the elliptic equation $-v''+ \gamma(v)= f'$ for a maximal monotone graph $\gamma$. Since, by accretivity, $\|\gamma(v)\|_{L^p} \leq \|f'\|_{L^p}$ for every $p$, and since $u'=v''+f'$ by the representation formula, one immediately concludes the estimate for $u'$ in terms of $f'$. The same kind of estimate also yields the  preservation of the BV norm in dimension one. 

In fact, it is easy to see that the approach used in Corollary \ref{estint} gives a  BV estimate in one dimension, as well.  
To this purpose, one multiplies the   equation \rife{div-form} by $(\vep+w)^{-\frac12}$ and  integrate by parts obtaining
\begin{align*}
& - \frac\de2\into \frac{|\nabla  w|^2}{(\vep+ w)^{\frac32}}dx
- \frac12 \into \frac{|\nabla  w|^2}{(\vep+ w)^{2}}dx
+ \frac12 \into \frac{\left(\nabla  w \nabla u\right)^2}{(\vep+ w)^{3}}dx
\\
& 
+ 2\lambda \into \frac w{(\vep+w)^{\frac12}} dx 
+2\de \, \into \frac{|D^2 u|^2}{{\sqrt{\vep +w}}}dx  +
2\into \frac{|D^2 u|^2}{{(\vep +w)}} 
\\
& \leq  \frac12 \into  \frac{|\nabla  w|^2}{(\vep+ w)^2}+ 2\into |\nabla f|\, \frac{|\nabla u|} {\sqrt{\vep +w}}dx\,.
\end{align*}
Rearranging terms and using $\frac{|\nabla u|} {\sqrt{\vep +w}}\leq 1$ we get
\be\label{rearr}
\begin{split}
& 2\de \into \left[\frac{|D^2 u|^2}{{\sqrt{\vep +w}}}- \frac14  \frac{|\nabla  w|^2}{(\vep+ w)^{\frac32}} \right]dx
+ \frac12 \into \frac{\left(\nabla  w \nabla u\right)^2}{(\vep+ w)^{3}}dx
 + 2\into \frac{|D^2 u|^2}{{(\vep +w)}} 
\\ & \quad 
+ 2\lambda \into \frac w{(\vep+w)^{\frac12}} dx   \leq   \into  \frac{|\nabla  w|^2}{(\vep+ w)^2}+ 2\into |\nabla f|\, dx\,.
\end{split}
\ee
Recall that $\nabla w= 2 D^2u \nabla u$. Hence the first integral is positive.
Moreover,  since in dimension one we have $|\nabla  w \nabla u| = |\nabla u|\, |\nabla w|$, by Young's inequality we estimate
$$
\frac{|\nabla  w|^2}{(\vep+ w)^2}   \leq  \,2\frac{|D^2 u|\, |\nabla u|\, |\nabla  w|}{(\vep+ w)^2}
\leq 
\, 2\frac{|D^2 u|^2}{(\vep+ w)} + \frac12\frac{|\nabla u \nabla  w|^2}{(\vep+ w)^3}
$$
Using the previous   inequality in \rife{rearr} we deduce
$$
2\lambda \into \frac w{(\vep+w)^{\frac12}} dx \leq 2\into |\nabla f|\, dx\, 
$$
which means that 
$$
\into \frac{|\nabla u_{\vep,\de}|^2}{\sqrt{\vep +|\nabla u_{\vep,\de}|^2} } \leq \frac1\lambda\, \|\nabla f\|_{L^1}\,.
$$
In particular, $u_{\vep,\de}$ satisfies
$$
\|\nabla u_{\vep,\de}\|_{L^1} \leq \frac1\lambda\, \|\nabla f\|_{L^1} + \sqrt \vep\, |\Omega|\,,
$$
which yields, when passing to the limit,  
\be\label{bv}
\| u\|_{BV(\Omega)} \leq   \| f\|_{BV(\Omega)}\,.  
\ee
We stress once more that this estimate is scale invariant, namely it does not depend on the possible coercivity parameter 
$\mu $ of the energy term in the functional \rife{TVmu}. Recall that minimizers of \rife{TVmu} also satisfy the standard estimate
$$
\into |Du|\, dx \leq \frac{\|f\|_{\elle2}^2}{2\mu}\,,
$$
but this latter one does actually depend on $\mu$, which is not the case for  \rife{bv}.
\end{rem}

\section{Proof of the regularity for the $BV$ minimizer.}

In this Section we conclude the proof of the results stated in the introduction, standing on the a priori estimates derived above.  We first deal with the global regularity results. 

\vskip0.5em

{\bf Proof of Theorem  \ref{lip} and  Theorem \ref{sob}. }

We consider $\de=\vep$ and we still call $u_\vep$ the unique minimizer of the functional
$$
J_\vep:= \into \left[\vep \frac{|\nabla u|^2}2+ \sqrt{\vep+ |\nabla u|^2}\right] dx +   \into \frac{|u-f|^2}2 dx \,. 
$$
Assume by now that $f\in W^{1,\infty}(\Omega)$. Then Corollary \ref{glob-lip} implies that
$$
\| \nabla u_{\vep}\|_\infty \leq  c_1 \, (\|\nabla f\|_\infty+ c_0)\,.
$$
Since by maximum principle we also have $\| u_{\vep}\|_\infty \leq   \|  f\|_\infty$, we deduce that $u_\vep$ is uniformly bounded in $W^{1,\infty}(\Omega)$. Thus, up to a subsequence, $u_\vep $ converges weakly to some $u\in W^{1,\infty}(\Omega)$ (and strongly in $\elle\infty$). By weak lower semicontinuity, we deduce the same estimate for $\nabla u$:
\be\label{stili}
\| \nabla u \|_\infty \leq    c_1 \, (\|\nabla f\|_\infty+ c_0)\,.
\ee
We are left to show that $u$ is the unique minimizer of \rife{TV}. To this purpose, let $v$ be any function in $BV(\Omega)\cap \elle2$; by standard properties of functions of bounded variation (see e.g. \cite{AFP}) there exists a sequence $v_n\in C^\infty(\Omega)$ such that 
$v_n\to v$ in $\elle 2$  and $\into |Dv_n|dx \to \into |Dv| \, dx$. Since $u_\vep$ is minimizer of $J_\vep$, we have
\begin{align*}
& \into \left[\vep \frac{|\nabla u_\vep|^2}2+ \sqrt{\vep+ |\nabla u_\vep|^2}\right] dx +   \frac12\into u_\vep^2 dx - \into f\, u_\vep \, dx
\\
& \qquad  \leq \into \left[\vep \frac{|\nabla v_n|^2}2+ \sqrt{\vep+ |\nabla v_n|^2}\right] dx + \frac12 \into v_n^2 dx - \into f\, v_n\, dx
\end{align*}
which implies
\be\label{preeps}
\begin{split}
& \into  |\nabla u_\vep| \,dx + \frac12 \into u_\vep^2 dx - \into f\, u_\vep \, dx
\\
& \qquad  \leq \into [\vep \frac{|\nabla v_n|^2}2+ \sqrt{\vep+ |\nabla v_n|^2}] dx + \frac12 \into v_n^2 dx - \into f\, v_n\, dx\,.
\end{split}
\ee
We let $\vep \to 0$ and we get
$$
\into  |\nabla u | \,dx + \frac12 \into u^2 dx - \into f\, u \, dx
  \leq \into |\nabla v_n|\,  dx + \frac12 \into v_n^2 dx - \into f\, v_n\, dx\,.
$$
Letting $n\to \infty$ we deduce that
$$
\into  |\nabla u | \,dx + \frac12 \into u^2 dx - \into f\, u \, dx
  \leq \into |Dv|\,  dx + \frac12 \into v^2 dx - \into f\, v \, dx\,.
  $$
This means that $u$ is the unique minimizer of \rife{TV}, and the whole sequence $u_\vep$ converges. Due to \rife{stili}, this proves Theorem \ref{lip}. The same proof applies for Theorem \ref{sob}: assuming $\Omega$ to be convex and $f\in W^{1,p}(\Omega)$, $p\geq 2$, we take some sequence of smooth functions $f_\vep $ converging to $f$ in $W^{1,p}(\Omega)$. By Corollary \ref{estint} we have
\be\label{presob}
\| \nabla u_{\vep}\|_{\elle p} \leq   \| \nabla f_\vep\|_{\elle p}\,.
\ee
A similar estimate $\|u_\vep\|_{\elle p}\leq    \|   f_\vep\|_{\elle p}$ holds   for $u_\vep$ as well. Thus $u_\vep$ is bounded in $W^{1,p}(\Omega)$ and, up to subsequences, converges weakly to some $u\in W^{1,p}(\Omega)$ 
and strongly in $\elle p$.   Since $p\geq 2$, there is no problem in repeating the same argument as above; passing to the limit first as $\vep \to 0$ and then as $n\to \infty$, we deduce that $u$ is the unique minimum of \rife{TV} and by lower semicontinuity $u$ preserves the estimate \rife{presob}. This proves Theorem \ref{sob}.
%
 \qed
  
  \vskip1em
The proof of Theorem \ref{liploc} slightly differs from the other theorems in two details. On one hand the local regularity assumed on $f$ does not imply any global bound on the approximating sequences; on another hand, the local estimate of Corollary \ref{est-loc} is obtained after letting the viscosity regularization vanish. So we will conclude by considering a standard convolution approximation of $f$ which initially enjoys global bounds and eventually   preserves the local properties of $f$.

\vskip1em
{\bf Proof of Theorem \ref{liploc}.}
Let $f_n:=\rho_n \star f$ be  a standard approximation of $f$ through convolution. Since  $f_n\in W^{1,\infty}(\Omega)$, we can apply to $f_n$ all the global bounds obtained before. In particular, let $u_{\vep,\de}$ be the unique solution to \rife{pbeps} corresponding to $f_n$; at fixed $n$, we know that $u_{\vep, \de}$ is bounded uniformly in $W^{1,\infty}(\Omega)$ and there are functions, called $u_\vep$ and $u_n$, such that 
 $u_{\vep, \de} \to u_\vep$ as $\de \to 0$ and $u_\vep \to u_n$ as $\vep \to 0$. The convergences are weak$-*$ in $W^{1,\infty}(\Omega)$ and strongly in $\elle\infty$. Using Corollary \ref{est-loc}, we deduce that
$$
\sup_{B_R(x_0)} \, |\nabla u_n| \leq   \left( \frac C{R^2} + \sup_{B_{\frac32 R}(x_0)} \, |\nabla f_n |\right)\,.
$$
By property of convolutions, if $f\in W^{1,\infty}(B_{2R}(x_0)) $ we have, for $n$ sufficiently large,
 $$
 \sup_{B_{\frac32 R}(x_0)} \, |\nabla f_n | \leq \sup_{B_{2 R}(x_0)} \, |\nabla f  |
$$
so we get
\be\label{loca}
\sup_{B_R(x_0)} \, |\nabla u_n| \leq   \left( \frac C{R^2} + \sup_{B_{2 R}(x_0)} \, |\nabla f |\right)\,.
\ee
We are left with the limit as $n\to \infty$. At fixed $n$, we know (as in the proof of Theorem \ref{lip}) that $u_n\in W^{1,\infty}(\Omega)$ and  is the unique minimizer of \rife{TV} corresponding to $f_n$. By minimizing property, $u_n$ is bounded in $BV(\Omega)$ and therefore relatively compact in $\elle 1$. Up to subsequences, we can assume that $u_n$ converges to $u\in BV(\Omega)\cap \elle2$ and by semicontinuity we conclude that $u$ is the unique minimizer of \rife{TV}. From \rife{loca}, it follows that $u\in W^{1,\infty}(B_R(x_0))$ and satisfies the same estimate.
\qed

\vskip1em
{\bf Acknowledgement.} The author warmly thanks H. Brezis for stimulating his interest in the problem.

\end{document}